\documentclass[reqno]{amsart}

\usepackage[all]{xy}

\numberwithin{equation}{section}

\theoremstyle{plain}
\newtheorem{lemma}{Lemma}[section]
\newtheorem{teo}[lemma]{Theorem}
\newtheorem{propo}[lemma]{Proposition}
\newtheorem{coro}[lemma]{Corollary}
\newtheorem{con}[lemma]{Conjecture}

\theoremstyle{definition}

\theoremstyle{remark}
\newtheorem{rem}[lemma]{Remark}

\newcommand{\pic}{{\rm Pic\thinspace}}
\newcommand{\bs}{{\rm Bs\thinspace}}
\newcommand{\p}{\mathbb{P}}

\newcommand{\z}{\mathbb{Z}}

\newcommand{\oc}{{\mathcal O}}
\newcommand{\ls}{{\mathcal L}}
\newcommand{\ms}{{\mathcal M}}
\newcommand{\ps}{{\mathcal P}}
\newcommand{\ci}{{\mathcal I}}

\newcommand{\rmap}{\dashrightarrow}

\newcommand{\vrk}{\hfill $\Box$}

\begin{document}

\title{Linear systems on generic $K3$ surfaces}
\author{Cindy De Volder}
\address{
Department of Pure Mathematics and Computeralgebra, Galglaan 2,
\newline B-9000 Ghent, Belgium} \email{cdv@cage.ugent.be}
\thanks{The first author is a Postdoctoral Fellow of
the Fund for Scientific Research-Flanders (Belgium)
(F.W.O.-Vlaanderen)}

\author{Antonio Laface}
\address{
Dipartimento di Matematica, Universit\`a degli Studi di Milano,
Via Saldini 50, \newline 20133 Milano, Italy }
\email{antonio.laface@unimi.it}
\thanks{The second author would like to thank the European Research and
Training Network EAGER for the support provided at Ghent
University. He also acknowledges the support of the MIUR of the
Italian Government in the framework of the National Research
Project ``Geometry in Algebraic Varieties'' (Cofin 2002)}
\keywords{Linear systems, fat points, generic $K3$ surfaces.}
\subjclass{14C20, 14J28.}
\begin{abstract}
In this paper we prove the equivalence of two conjectures on
linear systems through fat points on a generic $K3$ surface. The
first conjecture is exactly as Segre conjecture on the projective
plane. Whereas the second characterizes such linear system and can
be compared to the Gimigliano-Harbourne-Hirschowitz conjecture.
\end{abstract}
\maketitle

\section{Introduction}

In this paper we assume the ground field is algebraically closed
of characteristic 0.

With $S$ we always denote a smooth projective {\em generic} $K3$
surface, i.e. $\pic(S) \cong \z$.

Consider $r$ points in general position on $S$, to each one of
them associate a natural number $m_i$ called the {\em
multiplicity} of the point. Let $r_j$ be the number of $p_i$ with
multiplicity $m_i$ and let $H$ be the generator of $\pic S$.

For a linear system of curves in $|dH|$ with $r_j$ general base
points of multiplicity $m_j$ for $j=1\cdots k$, define its {\em
virtual dimension} $v$ as $\dim|dH|-\sum r_i m_i(m_i+1)/2$ and its
{\em expected dimension} by $e=\max\{v,-1\}$. If the dimension of
the linear system is $l$, then $ v\leq e\leq l$.

Observe that it is possible to have $e < l$, since the conditions
imposed by the points may be dependent. In this case we say that
the system is {\em special}.

Linear systems through general fat points on rational surfaces
have been studied by many authors (see e.g.
\cite{BS,AH,AG,BH,CM1,CM2}), but, as far as we know, no conjecture
concerning the structure of such systems on $K3$ surfaces has been
formulated.

Inspired by the article~\cite{CM1} by C.~Ciliberto and R.~Miranda,
we start with a Segre-like conjecture (conjecture~\ref{segre}),
and deduce conjecture~\ref{our}, which can be seen as a
translation of the Gimigliano-Harbourne-Hirschowitz conjecture
\cite[Conjecture 3.1]{CM1} to $K3$ surfaces.

In section~\ref{prel} we introduce some notation and definitions
and state the two conjectures. In the following section we prove
that the two conjectures are in fact equivalent. Finally, in
section~\ref{evi}, we prove some results which are in favor of
conjecture~\ref{segre}.

\section{Preliminaries}\label{prel}

Let $S$ be a generic $K3$ surface and let $H$ be the generator of
$\pic(S)$, then $H$ is ample, $H^2=2g-2 \geq 2$ and $h^0(H)=g+1$;
moreover $H$ is very ample if $g \geq 3$ and if $g=2$ $H$ defines
a double covering of $\p^2$ branched at an irreducible sextic (see
\cite[Proposition~3]{AM}).

Consider $p_1,\cdots, p_r$ points in general position on $S$, for
each one of these points fix a multiplicity $m_1,\cdots m_r$. With
$\ls = \ls^n(d,m_1,\cdots, m_r)$ we mean the linear system of
curves in $|dH|$ with multiplicity $m_i$ at $p_i$ for all $i$,
where $n = 2g-2$.

Let $Z=\sum m_ip_i$ be the $0$-dimensional scheme defined by the
multiple points and consider the exact sequence of sheaves:
\[
\xymatrix@1{ 0 \ar[r]  & \oc_S(dH) \otimes {\mathcal I}(Z) \ar[r]
& \oc_S(dH) \ar[r]  & \oc_Z \ar[r]  &  0 }
\]
where ${\mathcal I}(Z)$ is the ideal sheaf of $Z$. Taking
cohomology we obtain
\begin{equation}
\label{equx} v = h^0(\oc_S(dH) \otimes {\mathcal I}(Z)) -
h^1(\oc_S(dH) \otimes {\mathcal I}(Z)) -1,
\end{equation}
because $h^1(\oc_S(dH))=0$ (see e.g.~\cite{AM}).

Equation~(\ref{equx}) shows that a non-empty linear system $\ls$
is special if and only if $h^1(\ls) = 0$. We can now formulate the
following

\begin{con}\label{segre}
If $\ls$ on $S$ is non-empty and reduced, then it is non-special.
\end{con}

\begin{rem}\label{bertini}
Conjecture~\ref{segre} implies that a general element of a special
system is reduced, which is equivalent to saying that there exists
a curve $C$ such that $2C \subseteq \bs \ls$ (because of Bertini's
first theorem as stated in \cite[Theorem (4.1)]{SK}).
\end{rem}

\begin{con}\label{our}
Let $\ls$ and $S$ be as above, then
\begin{itemize}
\item[(i)]
    $\ls$ is special if and only if $\ls = \ls^4(d,2d)$ or $\ls = \ls^2(d,d^2)$ with $d \geq 2$;
\item[(ii)]
    if $\ls$ is non-empty then its general divisor has
    exactly the imposed multiplicities in the points $p_i$;
\item[(iii)]
    if $\ls$ is non-special and has a fixed irreducible component $C$
    then
    \begin{itemize}
    \item[a)] $\ls=\ls^2(m+1,m+1,m)=mC+\ls^2(1,1)$ with $C=\ls^2(1,1^2)$ or
    \item[b)] $\ls = 2C$ with $C \in \{ \ls^4(1,1^3),\allowbreak \ls^6(1,2,1),\allowbreak \ls^{10}(1,3)\}$ or
    \item[c)] $\ls = C$.
    \end{itemize}

\item[(iv)]
    if $\ls$ has no fixed components then either its general element is irreducible or
    $\ls=\ls^2(2,2)$.
\end{itemize}
\end{con}

\section{The equivalence of the two conjectures}

It is clear that conjecture~\ref{our} implies
conjecture~\ref{segre}, and we will now show that actually they
are equivalent.

For the rest of this section we assume that conjecture~\ref{segre}
is true.

If $\ls$ and $S$ are as above, let $S'$ denote the blowing-up of
$S$ along the points $p_1,\ldots,p_r$, and let $E_i$ be the
exceptional divisor on $S'$ corresponding to $p_i$. Then the
canonical class $K_{S'}$ of $S'$ is equal to $\sum_{i=1}^r E_i$
and $\pic S'$ is generated by $\{H,E_1,\ldots,E_r\}$, where, by
abuse of notation, $H$ also denotes the pullback of $H$ on $S'$.

Let $D$ be a divisor on $S'$, such that $|D| = |tH - \sum_{i=1}^r
l_i E_i|$ with $t \geq 0$ for all $i$. Then $\chi(D) =
\frac{1}{2}(D^2 - D K_{S'}) +2$, and we define the virtual
dimension of $D$ as
$$v(D) = \chi(D) - h^2(D) -1 = h^0(D) - h^1(D) - 1. $$
By duality, $h^2(D) = 0$, unless $|D| = |\sum_{i \in I} E_i|$,
with $I \subseteq \{1,\ldots,r\}$. Also note that if $t >0$ then
$|D|$ corresponds to the system $\ls^n(t,l_1,\ldots,l_r)$ on $S$
and $v(D) = v(\ls^n(t,l_1,\ldots,l_r))$ (see
equation~(\ref{equx})). By abuse of notation we then also denote
$|D|$ by $\ls^n(t,l_1,\ldots ,l_r)$. Moreover, if $C$ and $C'$ are
two curves on $S$, then by $C C'$ we mean the intersection
multiplicity of their strict transforms on $S'$.

\begin{lemma}\label{irr}
Let $\ms$ be a linear system on $S$ without fixed components then
either its general element is irreducible or $\ms=\ls^2(2,2)$.
\end{lemma}

\begin{proof}
Because of Bertini second theorem as stated in~\cite[5.3]{SK} the
general element of $\ms$ is reducible if and only it is composite
with a pencil $\ps$. Let $M$ and $P$ be the strict transforms on
$S'$ of general elements of $\ms$ and $\ps$ respectively. This
means that $|M|=|lP|$ with $l=\dim|M|\geq 2$. By
remark~\ref{bertini} we can say that $|M|$ and $|P|$ are
non-special which gives $v(M)=l$ and $v(P)=1$. The second equality
implies that $PK_{S'}=P^2$ so the first is equivalent to $P^2=2/l$
which gives $l=2$ and $P^2=PK_{S'}=1$. This means that $P=dH-E_i$
and $1=P^2=nd^2-1$ which is only possible if $n=1$ and $d=1$, i.e.
$\ms=\ls^2(2,2)$.
\end{proof}

Observe that for any divisors $A,B$ on $S'$ we have
$\chi(A+B)=\chi(A)+\chi(B)+AB-2$, hence if $h^2(A)=h^2(B)=0$ then
\begin{equation}\label{add}
v(A+B)=v(A)+v(B)+AB-1.
\end{equation}

\begin{lemma}\label{exc}
Let $|D|$ be the linear system on $S'$ corresponding to a linear
system $\ls = \ls^n(d,m_1,\ldots,m_r)$ on $S$. Then $E_i
\not\subseteq \bs |D|$ for all $i = 1,\ldots,r$.
\end{lemma}

\begin{proof}
Assume the statement is false, then we can write $|D| = \sum l_i
E_i + F + |D'|$, with $E_i \not\subseteq \bs |D'|$, $F$ the strict
transform of the fixed components of $\ls$ and $|D'|$ without
fixed components.
\\
Assume there exists an $i$ such that $l_i >0$, then $E_i (D' + F)
>0$ (otherwise $m_i = D E_i = -l_i <0$).
\\
Assume that $E_i F > 0$. Let $C \subseteq F$ be an irreducible
divisor with $C E_i >0$. By conjecture~\ref{segre} $\dim C = v(C)
= 0$, and $C+E_i$ still has dimension 0, because it is contained
in $\bs |D|$. On the other hand, $C+E_i$ is non-special (since $C$
is non-special), so $v(C+E_i) = v(C)$. This implies $\chi (C) =
\chi (C+E_i)$ which is equivalent to $C E_i = 0$ and thus
contradicts our assumption.
\\
So we get that $E_i D' >0$. By conjecture~\ref{segre} $\dim |D'| =
v(D')$, and, as before, $|D'| +E_i$ is non-special. So $D' E_i =
0$ which again contradicts our assumption.
\end{proof}

\begin{lemma}\label{table}
If $C$ is an irreducible divisor on $S'$ such that $v(C)=0$ and
$C^2 \leq 1$, then $C$ is one of the following:

\begin{center}
\begin{tabular}{|c|ccc|}
\hline
  & \multicolumn{3}{c|}{\vspace{-3mm}} \\
$C^2$ & \multicolumn{3}{c|}{$C$} \\
\hline
  & \multicolumn{3}{c|}{\vspace{-3mm}} \\
$\leq -1$ & \multicolumn{3}{c|}{$\emptyset$} \\
$0$ & $\ls^4(1,2)$      & $\ls^2(1,1^2)$        & \\
$1$ & $\ls^4(1,1^3)$    & $ \ls^6(1,2,1)$   & $\ls^{10}(1,3)$ \\
\hline
\end{tabular}
\end{center}

\end{lemma}

\begin{proof}
If $h^2(C) > 0$, then $C = E_i$ for some $i$, in which case $v(C)
= 0$ and $C^2 = - 1$.
\\
If $h^2(C) = 0$, then $v(C) =0$ implies that $C K_{S'} = C^2 + 2$.
But $C^2 + 2 = p_a(C) \geq 0$, so $C^2 \geq -2$.
\\
In case $C^2 \leq -1$, $C K_{S'} \leq 1$ so $C = tH - aE_i$ with
$a \in \{0,1\}$, which gives $t^2n - a \leq -1$ and this is not
possible if $t >0$.
\\
In case $C^2 = 0$ then $C K_{S'} = 2$, so either $C = tH - 2E_i$
or $C = tH - E_i - E_j$. In the first case, $C^2 = nt^2 - 4 =0$,
which is only possible if $n=4$ and $t=1$. In the latter case,
$C^2 = nt^2 -2 =0$, which is only possible if $n=2$ and $t=1$.
\\
In case $C^2 = 1$ then $C K_{S'} = 3$, so $C$ is of the following
types $tH- 3E_i$, $tH - 2E_i - E_j$ or $tH - E_i -E_j -E_k$. In
the first case, $C^2 = nt^2 - 9 = 1$, which is only possible if
$n=10$ and $t=1$. In the second case, $C^2 = nt^2 - 5 = 1$, which
is only possible if $n=6$ and $t=1$. And in the latter case, $C^2
= nt^2 - 3 = 1$, which is only possible if $n=4$ and $t=1$.
\end{proof}

\begin{propo}\label{fixed}
Let $\ls=\ls^n(d,m_1,\ldots,m_r)$ be a linear system on $S$ and
assume that there exist distinct irreducible curves $C_i$ and
$D_j$ such that the fixed part of $\ls$ is given by
$\sum_{i=1}^a\mu_i C_i + \sum_{i=1}^bD_i$, where $\mu_i\geq 2$,
then either $\dim \ls = 0$ and it is one of the following:

\begin{center}
\begin{tabular}{lccl}
 $m C$                      & $C^2=0$       & $v(C)=0$      & special for $m\geq 2$ with $v(\ls)= 1-m$ \\
 $2C$                               & $C^2=1$       & $v(C)=0$      & non-special \\
 $D$                                &   $D^2\geq 0$ &   $v(D)=0$        & non-special \\
\end{tabular}
\end{center}

or $\ls=\ls^2(m+1,m+1,m)=m\ls^2(1,1^2)+\ls^2(1,1)$ which is
non-special. Note that in the latter case $m\ls^2(1,1^2)$ is the
fixed part of $\ls$ and $\ls^2(1,1)$ is its free part of dimension
$1$.

\end{propo}

To simplify the proof of this proposition, we first give two
lemmas.

\begin{lemma}\label{int}
With the same situation of proposition~\ref{fixed} we have
$C_iC_j=C_iD_j=D_iD_j=1$ and $C_i^2\leq 1$.
\end{lemma}

\begin{proof}
Since $C_i+C_j\subset\bs(\ls)$ then $\dim|C_i+C_j|=0$ and
conjecture~\ref{segre} implies that $v(C_i+C_j)=0$. Using the same
argument one also obtains that $v(C_i)=v(C_j)=0$, so
equation~(\ref{add}) implies that $C_iC_j=1$. In the same way one
can proof that $C_iD_j=1$ and $D_iD_j=1$.

From $v(C_i)=0$ we obtain $C_iK_{S'}=C_i^2+2$; since
$2C_i\subset\bs(\ls)$ this implies that $v(2C_i)\leq 0$ which
gives $C_i^2\leq 1$.
\end{proof}

\begin{lemma}\label{pair}
Let $A$ and $B$ be the strict transforms on $S'$ of two distinct
irreducible curves on $S$ then one of the following holds:
\begin{itemize}
\item[] $A=\ls^2(1,1^2)$ and $B$ is an irreducible element of $\ls^2(1,1)$; or
\item[] $AB\neq 1$.
\end{itemize}
\end{lemma}

\begin{proof}
Consider the following exact sequence:

\[
\xymatrix@1{ 0 \ar[r] &   \oc_{S'}(A-B) \ar[r] &  \oc_{S'}(A)
\ar[r] &  \oc_B(A) \ar[r] &  0. }
\]

Using conjecture~\ref{segre}, we obtain $h^1(\oc_{S'}(A))=0$ and
as $h^2(\oc_{S'}(A))=0$ the preceding sequence implies that
$h^2(\oc_{S'}(A-B))=h^1(\oc_{B}(A))$. Moreover, by Riemann-Roch
and $h^0(\oc_{B}(A))=1$ we obtain that $h^1(\oc_{B}(A))=p_a(B)-1$.

Observe that $p_a(B)=(B^2+BK_{S'})/2+1\geq 2$; indeed, $B^2\geq 0$
and $BK_{S'}\geq 0$ but one can immediately verify that they can
not be both $0$.

The preceding calculation shows that, by Serre duality:
\[
\dim|K_{S'}+B-A|\geq p_a(B)-2.
\]
By interchanging the roles of $A$ and $B$ we obtain also
$\dim|K_{S'}+A-B|\geq p_a(A)-2$. This implies that
$\dim|2K_{S'}|\geq (p_a(A)-2) + (p_a(B)-2)|$. Since
$\dim|2K_{S'}|=0$ this means that $p_a(A)=p_a(B)=2$.

From the exact sequence we also see that
\[
h^0(\oc_{S'}(A))\leq h^0(\oc_{S'}(A-B))+h^0(\oc_{B}(A));
\]
on the other hand if $h^0(\oc_{S'}(A-B))\geq 2$ then
$\dim|K_{S'}|\geq\dim|K_{S'}+B-A|+\dim|A-B|\geq 1$ which is not
possible. So this implies that $\dim|A|\leq 1$.

We already know that $p_a(A)=2$ and the preceding analysis shows
that $0\leq v(A)\leq 1$. In case $v(A)=0$ we obtain $A^2=0$ so by
lemma~\ref{table} $A=\ls^4(1,2)$ or $A=\ls^2(1,1^2)$. Otherwise
$v(A)=1$ and $A^2=1$, this gives $AK_{S'}=1$ hence the only
possibility is given by $A=\ls^2(1,1)$.

The same considerations are also true for $B$, hence the only
possible pair $A,B$ such that $AB=1$ is given by: $A=\ls^2(1,1^2)$
and $B$ is an irreducible element of $\ls^2(1,1)$.
\end{proof}

\begin{proof}[Proof of proposition~\ref{fixed}]

Let $|L|$ be the linear system on $S'$ with $L$ the strict
transform of a general element of $\ls$ on $S$. Then
\[
|L|=\sum_{i=1}^a\mu_i C_i + \sum_{i=1}^b D_i + |M|,
\]
where $|M|$ is without fixed components and $C_i$ and $D_j$ are
distinct irreducible curves. And, by lemma~\ref{exc},
$\sum_{i=1}^a\mu_i C_i + \sum_{i=1}^bD_i$ is the strict transform
of the fixed part of $\ls$.
\\
Because of lemmas~\ref{int} and \ref{pair} $a+b = 1$.
\\
Assume that $|M|$ is non-trivial, i.e. $\dim \ls >0$. Let $C$ be
an irreducible fixed component of $|L|$, then, by
conjecture~\ref{segre}, $v(M) = v(M + C)$, which implies that $M C
=1$. By lemma~\ref{irr} either the general element of $|M|$ is
irreducible or $|M| = \ls^2(2,2)$.
\\
In the latter case a general element of $|M|$ can be written as
$M_1 + M_2$ with $M_i \in \ls^2(1,1)$, so $M C$ cannot be equal to
1 as it is always even.
\\
In the first case, by lemma~\ref{pair}, the only possibility is $C
= \ls^2(1,1^2)$ and $|M| = \ls^2(1,1)$, i.e. $\ls =
\ls^2(m+1,m+1,m)=m\ls^2(1,1^2)+\ls^2(1,1)$. In order to see that
the last equality is true, we just have to note that $\dim \ls =
1$. Indeed, by specializing the two general points of $S$ to
points on the ramification divisor, the obtained system
corresponds to $\oc_{\p^2}(m+1) \otimes \ci(Z)$ with $Z = (m+1)p_1
+ m p_2$, which, obviously has dimension 1.
\\
Now assume that $\dim \ls =0$. In case $b=1$, by
conjecture~\ref{segre}, $v(D) =0$, and, by lemma~\ref{table}, we
know that $D^2 \geq 0$. In case $a=1$, by lemma~\ref{int} we know
that $\ls = m C$ with $C^2 \leq 1$ and $v(C) =0$ (because of
conjecture~\ref{segre}). If $C^2 = 0$ then $CK_{S'} = 2$, so
$v(mC) = 1 - m$. If $C^2 = 1$ then $CK_{S'} = 3$, so $v(mC) =
m(m-3)/2 + 1$ thus $v(mC) \leq 0$ implies $m =2$.
\end{proof}

Note that proposition~\ref{fixed} and lemma~\ref{table} imply (i)
and (iii) of conjecture~\ref{our}, part (ii) follows from
lemma~\ref{exc} and lemma~\ref{irr} implies (iv); so we proved the
following

\begin{teo}
Conjecture~\ref{segre} implies conjecture~\ref{our}.
\end{teo}

\section{Results in favor of conjecture~\ref{segre}}\label{evi}

In this section we will list some results which leads us to
believe that conjecture~\ref{segre} is true.

\begin{teo}\label{node}
Let $\ls$ be a non-special linear system on a smooth projective
surface $X$ such that $\ls \otimes \ci(2p)$ is special for a
general point $p \in X$, then $\ls \otimes \ci(2p)$ has a double
fixed component through $p$.
\end{teo}

\begin{proof}
We may assume that $\ls$ has no fixed components, because
otherwise we can consider $\ls - F$, where $F$ is the fixed
divisor of $\ls$. Let $n = \dim \ls$ and consider the rational map
$\varphi: X \rmap \p^n$ corresponding to $\ls$. Saying that $\ls
\otimes \ci(2p)$ is special for a general point $p \in X$, means
that the image $X' = \varphi(X)$ has to be a curve. Indeed, the
speciality implies that an hyperplane which contains $p' =
\varphi(p)$ and one tangent direction $\tau \in T_{p'}(X')$ has to
contain the whole tangent space $T_{p'}(X')$. Because $\ls$ is
given by $\varphi^*(\oc_{\p^n}(1) |_{X'})$, the general divisor of
$\ls$ can be written as $\varphi^*(\sum n_i p_i) = \sum n_i F_i$,
with $\sum n_i = \deg X'$. So any divisor of $\ls \otimes \ci(2p)$
contains $2 F_p$ with $F_p = \varphi^*(\varphi(p))$, i.e. $2 F_p
\subset \bs (\ls \otimes \ci(2p))$.
\end{proof}

If it is possible to find reduced curves $C_1$, $C_2$ on a generic
$K3$ surface such that $C_1$ and $C_2$ have no common components,
$v(C_1) = v(C_2) = 0$ and $v(C_1 + C_2) <0$, then this would imply
that conjecture~\ref{segre} is false because $C_1 + C_2$ would be
a special system with no multiple components. We will show that
such curves can not exist in the following

\begin{propo}\label{C+C'}
Let $S$ be a generic $K3$ surface with $H^2 = n$, let $C$ and $C'$
be curves on $S$ with $v(C) = v(C') =0$, then $v(C + C') \geq 0$
unless $C = C' \in \{ \ls^2(1,1^2), \ls^4(1,2) \}$.
\end{propo}

\begin{proof}
Let $C \in \ls^{n}(d,m_1,\ldots,m_r)$ and $C' \in
\ls^{n}(d',m'_1,\ldots,m'_r)$. Then $v(C) = 0$, resp. $v(C')=0$,
implies that $n d^2 = \sum_{i=1}^r m_i(m_i+1) -2$, resp. $n {d'}^2
= \sum_{i=1}^r m'_i(m'_i+1) -2$, so
$$(n d d')^2 = (\sum_{i=1}^r m_i(m_i+1) -2)(\sum_{i=1}^r m'_i(m'_i+1) -2).$$
Moreover $v(C + C') \geq 0$ is equivalent to $n d d' >
\sum_{i=1}^r m_i m'_i$, which gives
\begin{equation}\label{int ineq}
(\sum_{i=1}^r m_i(m_i+1) -2)(\sum_{i=1}^r m'_i(m'_i+1) -2) >
(\sum_{i=1}^r m_i m'_i)^2.
\end{equation}
Note that, by Schwartz inequality,
$$(\sum m^2_i)(\sum {m'}^2_i) - (\sum m_i m'_i)^2 \geq 0,$$
and we denote the left hand side by $t$. To simplify notation we
also denote $a= \sum m^2_i$, $a'= \sum {m'}^2_i$, $b= \sum m_i$
and $b' = \sum m'_i$. Then the inequality~(\ref{int ineq}) is
equivalent to
$$
t + (a - 2) (b + b' -2) + a' (b - 2) > 0.
$$
Because $v(C)= v(C') =0$, we know that $a, b \geq 2$, so the
preceding inequality is true unless $a = b = 2$ and $t = 0$. We
know that $t = 0$ if and only if there exists a constant $c$ such
that $m'_i = c m_i$ for all $i$. On the other hand $a = b = 2$
implies that $C, C' \in \{ \ls^n(d,1^2) , \ls^n(d,2)\}$, and $0  =
v(\ls^n(d,1^2)) = nd^2/2 - 1$ is only possible for $n=2$ and
$d=1$, while $0  = v(\ls^n(d,2)) = nd^2/2 - 2$ is only possible
for $n=4$ and $d=1$.
\end{proof}

\begin{coro}
Lemma~\ref{int} is true without assuming conjecture~\ref{segre}.
\vrk
\end{coro}

\begin{coro}
With the same assumptions as in proposition~\ref{fixed} we have
that $a+b \leq 2$ (without assuming conjecture~\ref{segre}).
\end{coro}

\begin{proof}
Assume that $a + b \geq 3$ and let $R_1, R_2, R_3 \subseteq \bs
\ls$ be three distinct irreducible curves. Then $v(R_i)=0$ and,
since $R_i + R_j \subseteq \bs \ls$, proposition~\ref{C+C'}
implies that $v(R_i + R_j) = 0$; which in turn is equivalent to
$R_i R_j =1$. But, in the same way we can prove that $R_1 (R_2 +
R_3) = 1$, which contradicts $R_i R_j =1$.
\end{proof}

\begin{rem}
Observe that lemma~\ref{table} is true without assuming
conjecture~\ref{segre}; and the same is true for lemma~\ref{irr},
with just minor changes to the proof.
\end{rem}

\bibliographystyle{alpha}

\begin{thebibliography}{May72}

\bibitem[CM98]{CM2}
Ciro Ciliberto and Rick Miranda.
\newblock Degenerations of planar linear systems.
\newblock {\em J. Reine Angew. Math.}, 501:191--220, 1998.

\bibitem[CM01]{CM1}
C.~Ciliberto and R.~Miranda.
\newblock The {S}egre and {H}arbourne-{H}irschowitz conjectures.
\newblock In {\em Applications of algebraic geometry to coding theory, physics
  and computation (Eilat, 2001)}, volume~36 of {\em NATO Sci. Ser. II Math.
  Phys. Chem.}, pages 37--51. Kluwer Acad. Publ., Dordrecht, 2001.

\bibitem[Gim89]{AG}
Alessandro Gimigliano.
\newblock Regularity of linear systems of plane curves.
\newblock {\em J. Algebra}, 124(2):447--460, 1989.

\bibitem[Har86]{BH}
Brian Harbourne.
\newblock The geometry of rational surfaces and {H}ilbert functions of points
  in the plane.
\newblock In {\em Proceedings of the 1984 Vancouver conference in algebraic
  geometry}, volume~6 of {\em CMS Conf. Proc.}, pages 95--111, Providence, RI,
  1986. Amer. Math. Soc.

\bibitem[Hir89]{AH}
Andr{\'e} Hirschowitz.
\newblock Une conjecture pour la cohomologie des diviseurs sur les surfaces
  rationnelles g\'en\'eriques.
\newblock {\em J. Reine Angew. Math.}, 397:208--213, 1989.

\bibitem[Kle98]{SK}
Steven~L. Kleiman.
\newblock Bertini and his two fundamental theorems.
\newblock {\em Rend. Circ. Mat. Palermo (2) Suppl.}, (55):9--37, 1998.
\newblock Studies in the history of modern mathematics, III.

\bibitem[May72]{AM}
Alan~L. Mayer.
\newblock Families of $K3$ surfaces.
\newblock {\em Nagoya Math. J.}, 48:1--17, 1972.

\bibitem[Seg62]{BS}
Beniamino Segre.
\newblock Alcune questioni su insiemi finiti di punti in geometria algebrica.
\newblock In {\em Atti Convegno Internaz. Geometria Algebrica (Torino, 1961)},
  pages 15--33. Rattero, Turin, 1962.

\end{thebibliography}

\end{document}